\documentclass{amsart}
\def\platz{\par\medbreak}

   {\end{trivlist}\platz}

\def\qed{\vbox{\hrule
  \hbox{\vrule\hbox to 5pt{\vbox to 8pt{\vfil}\hfil}\vrule}\hrule}}
\def\endproof{\hfill\qed\medbreak}

\newtheorem{theorem}{Theorem}%[section]
\newtheorem{proposition}[theorem]{Proposition}

\newtheorem{corollary}[theorem]{Corollary}

\newtheorem{example}[theorem]{Example}

\newcommand\calS{{\mathcal S}}

\newcommand\calX{{\mathcal X}}

\newcommand {\R}{\mathbb{R}}

\newcommand {\grad}{\nabla}

\usepackage{amsrefs}

\newcommand{\hclass}{h_{\text{class}}}
\renewcommand{\div}{\operatorname{div}}
\newcommand{\LinfOmega}{L^\infty(\Omega,\R^n)}
\newcommand{\Xdiv}{\calX_{\div}}
\newcommand{\XBV}{\calX_{\operatorname{BV}}}
\newcommand{\phit}{{\phi^{-1}(t)}}
\newcommand{\philt}{{\{\phi<t\}}}
\newcommand{\cost}{\operatorname{cap}}
\renewcommand{\value}{\operatorname{value}}
\newcommand{\fmax}{f_{\text{max}}}\newcommand{\Smin}{S_{\text{min}}}
\newcommand{\rhotilde}{\tilde{\rho}}

\begin{document}

\title[First eigenvalue and Max Flow Min Cut]{The first eigenvalue of the Laplacian, isoperimetric constants, and the Max Flow Min Cut
Theorem}
\author{Daniel Grieser}
\address{Mathematisches Institut\\ Universit\"at Bonn\\ Beringstr.\ 1\\ 53115 Bonn\\
Germany}
\email{grieser@math.uni-bonn.de}
%\urladdr{}
\date{\today}
\begin{abstract}
We show how 'test' vector fields may be used to give lower bounds for the Cheeger constant of a Euclidean domain (or Riemannian manifold with boundary), and hence for the lowest eigenvalue of the Dirichlet Laplacian on the domain. Also, we show that a continuous version of the classical Max Flow Min Cut Theorem for networks implies that Cheeger's constant may be obtained precisely from such vector fields. Finally, we apply these ideas to reprove a known lower bound for Cheeger's constant in terms of the inradius of a plane domain.
\end{abstract}

\maketitle

%Keywords also: Osserman/Hayman inequality, Optimization

\section{Introduction}
For a domain $\Omega\subset \R^{n}$  the \emph{fundamental frequency} is defined by
\begin{equation}
\label{fundfreq}
\lambda_{\Omega} = \inf_{u\in C^\infty_{0}(\Omega)} R(u),\quad R(u)=
                 \frac{\int_{\Omega}|\nabla u|^2}{\int_{\Omega}u^2}.
\end{equation}
If $\Omega$ is bounded and has a Lipschitz boundary then this is the smallest eigenvalue of
the Laplacian $-\sum_i{\partial^2}/\partial x_i^{2}$ on $\Omega$ with Dirichlet boundary conditions. The minimum is attained by the corresponding eigenfunction, which lies in
$H^1_{0}(\Omega)\cap C^\infty(\Omega)$.

For most domains it is impossible to determine $\lambda_{\Omega}$ precisely, so it
is a fundamental problem to give estimates in terms of the geometry of $\Omega$. Upper estimates can be obtained by choice of any 'test function' $u$ in \eqref{fundfreq}.

It is less clear how to obtain lower estimates for $\lambda_{\Omega}$. One such
estimate was given by Cheeger in \cite{Che:LBSEL}: Define the \emph{Cheeger constant} by
\begin{equation}
\label{Cheegerconst}
h_{\Omega} = \inf_{S\subset\Omega} \frac {|\partial S|} {|S|}.
\end{equation}
The infimum is taken over open subsets $S$.
The absolute value signs denote $(n-1)$-dimensional Hausdorff measure in the
numerator and $n$-dimensional Hausdorff (= Lebesgue) measure in the denominator.
Then \emph{Cheeger's inequality}
says that
\begin{equation}
\label{Cheegerineq}
\lambda_{\Omega} \geq h_{\Omega}^2/4.
\end{equation}

The Cheeger constant is sometimes called an isoperimetric constant since it resembles the classical isoperimetric constant
$\hclass=\inf_{S\subset\Omega}\frac{|\partial S|^{n/(n-1)}}{|S|}$.
$\hclass$ is scale invariant and in fact independent of $\Omega$, and
the solution of the classical isoperimetric problem is that the mimimizers are
precisely the balls.

In contrast, $h_\Omega$ depends strongly on $\Omega$ and clearly scales as
$h_{r\Omega}=r^{-1}h_\Omega$. In general it is difficult to determine $h_\Omega$
precisely. It is known that minimizers $S$ exist if $\partial\Omega$ is Lipschitz, and that
$\partial S\cap\Omega$ is smooth (if $n\leq7$) and has constant mean curvature
(see \cite{FriKaw:IEFEPLOCC}, Theorem 8 and Remark 9). This may be used to determine $h_\Omega$ explicitly in some cases, for example for polygons (see \cite{KawLac:CCSCSP}),
and to give conditions when $\Omega$ is itself a minimizer (see \cite{AltCasCha:CCCSR}, \cite{Giu:ESPMCEUWBC};
such sets are called calibrable in the image processing literature).
In general, one may hope for estimates on $h_\Omega$ in terms of geometric data.

Again, upper estimates for $h_\Omega$ are obtained by using a suitable 'test domain' $S$, while it is less obvious how to obtain lower estimates.

One purpose of the present note is to point out a very simple idea how to obtain a lower estimate for $h_{\Omega}$:

\begin{proposition}\label{prop}
Let $V:\Omega\to\R^n$ be a smooth vector field on $\Omega$, $h\in\R$, and assume
\begin{eqnarray}
|V| &\leq & 1  \label{Vlength}\\
\div V&\geq & h \label{Vdiv},
\end{eqnarray}
both pointwise in $\Omega$.

Then $h_\Omega\geq h$.
\end{proposition}
\proof
Clearly, one may restrict to sets $S$ with smooth boundary in \eqref{Cheegerconst}.
For such $S$ we have, by Green's formula and \eqref{Vlength}, \eqref{Vdiv},
\begin{equation}\label{dSVS}
|\partial S|\geq \int_{\partial S} V\cdot dn = \int_S \div V \geq h|S|.
\end{equation}
\endproof

Proposition \ref{prop} seems to be little known in the geometric analysis community, although it is implicit
in McKean's proof of lower bounds for $\lambda_\Omega$ in case $\Omega$ is a complete, simply connected Riemannian manifold of strictly negative curvature
\cite{McK:UBSDMNC} (here, $V$ is taken as gradient of the distance to a fixed point). See \cite{CheLeu:EESWBMCHS} and \cite{BesMon:EESWLBMC} for other applications of the same idea, and also
\cite{BesMon:EBTGA} for another lower bound on $\lambda_\Omega$ in terms of vector fields.

\begin{example} \label{ex}
Let $\Omega=\{x\in\R^n:\, |x|<1\}$. Since $|\partial\Omega|/|\Omega| = n$, we have
$h_{\Omega}\leq n$. The vector field $V(x)=x$ has $|V|\leq 1$ and $\div V\equiv n$, so that, in fact, $h_{\Omega}= n$.
\end{example}

It is important for the sequel to allow non-smooth vector fields.
Consider the following classes:
\begin{eqnarray*}
\Xdiv (\Omega) & = & \{V\in\LinfOmega:\, \div V \in L^2(\Omega)\}  \\
\XBV(\Omega)  & = & \{V\in\LinfOmega:\, V \text{ has bounded variation}\}.
\end{eqnarray*}

$\div V$ is understood in the sense of distributions. Recall that, by definition, $V\in\LinfOmega$ has \emph{bounded variation} if all of its first derivatives $\partial V_i/\partial x_j$ (in the sense of distributions) are (signed) measures.
For a vector field  $V\in \LinfOmega$, \eqref{Vlength} is meant to hold almost everywhere and \eqref{Vdiv} in the sense of distributions. If $V\in\Xdiv\cup\XBV$ then $\div V$ is a measure, so \eqref{Vdiv} then holds also in the sense of measures.

Below, the class $\Xdiv$ will occur in the context of the Max Flow Min Cut Theorem,
and the class $\XBV$ will arise for vector fields defined via the distance function.
\medskip

\noindent{\bf Addendum to Proposition \ref{prop}}. \emph{Proposition \ref{prop} holds for vector fields $V\in \Xdiv\cup\XBV$.}

\proof
The proof \eqref{dSVS} still works, since for such $V$ one may define a 'restriction to the boundary'  $V_{|\partial S}$ (for open $S\subset\Omega$ with Lipschitz boundary), which satisfies $\|V_{|\partial S}\|_{L^\infty(\partial S)} \leq \|V\|_{L^\infty(S)}$ and
Green's formula.
For  $V\in\Xdiv$ this is shown in  \cite{Anz:PBMBFCC}.
For $V\in\XBV$ this follows from results in  \cite{EvaGar:MTFPF}, Section 5.3. Theorem 1 there states that a function $f\in L^1(S)$ of bounded variation on a Lipschitz domain $S\subset\R^n$ has a well-defined restriction to the boundary $f_{|\partial S}\in L^1(\partial S)$ satisfying, for any $W\in C^1(\R^n,\R^n)$,
$ \int_{\partial S} f_{|\partial S} W\cdot dn = \int_S f\div W + \int_S W\cdot\grad f.$
Applying this to $f=V_i$ and $W\equiv e_i$ (the $i$th standard unit vector)
for each $i=1,\dots,n$ and summing over $i$ yields Green's formula for $V$.
Also, by Theorem 2 loc.cit.\ one has $\|V_{|\partial S}\|_{L^\infty(\partial S)} \leq \|V\|_{L^\infty(S)}$.
\endproof

It is a remarkable fact that the estimate in Proposition \ref{prop} is sharp. That is, $h_\Omega$ may be characterized using
vector fields:

\begin{theorem}\label{MFMCcont}
We have
$$h_\Omega = \sup \{h:\, \exists V \text{ satisfying } \eqref{Vlength},\eqref{Vdiv}\},$$
where the supremum is taken over smooth vector fields $V$ on $\Omega$.

If $\partial\Omega$ is Lipschitz then there is a maximizing $V\in\Xdiv$.
\end{theorem}

Theorem \ref{MFMCcont} may be regarded as a continuous version of the classical Max Flow Min Cut Theorem for networks.
It was first proved by Strang \cite{Str:MFTD} in two dimensions and by Nozawa \cite{Noz:MFMCTAN} in general (see Theorem 4.4 there; Nozawa actually establishes a maximizing $V$ with $\div V\in L^n(\Omega)$). We explain the relation to the Max Flow Min Cut Theorem and  sketch the proof of Theorem \ref{MFMCcont} in Section \ref{secMFMC}.

Given the Theorem, one can prove Cheeger's inequality easily:

\proof (of Cheeger's inequality \eqref{Cheegerineq} using Theorem \ref{MFMCcont}.)
If $u\in C_0^\infty(\Omega)$ and $V$ is a smooth vector field satisfying
\eqref{Vlength}, \eqref{Vdiv} then, using Green's formula,
(with all integrals over $\Omega$)
\begin{multline*}
h\int u^2 \leq \int(\div V)\,u^2 = -\int V\cdot \grad(u^2) \\
          \leq 2\int |u|\,|\grad u| \leq 2\sqrt{\int u^2}\sqrt{\int|\grad u|^2},
\end{multline*}
so $R(u)\geq h^2/4$.
Taking a sequence $V_k$ with $h_k=\inf\div V_k$ approaching $h_\Omega$ one obtains $R(u)\geq h_\Omega^2/4$ for all $u\in C_0^\infty(\Omega)$ and therefore
\eqref{Cheegerineq}.
\endproof

This is not a substantially new proof of Cheeger's inequality: Cheeger's original proof is essentially a similar estimate, plus a clever use of the coarea formula applied to $u^2$. But the proof of Theorem \ref{MFMCcont} also relies on the coarea formula (see Section \ref{secMFMC})!

We remark that Proposition \ref{prop} and Theorem \ref{MFMCcont} extend directly to Riemannian manifolds with boundary, although (for Theorem \ref{MFMCcont}) this is not stated explicitly
in \cite{Noz:MFMCTAN}.
The relationship of Cheeger's inequality and Max Flow Min Cut Theorems was first noted by Alon \cite{Alo:EE} in the context of graphs (see also \cite{Chu:SGT}).

In this paper we consider Cheeger's inequality for the Dirichlet problem only since the case of closed manifolds or the Neumann problem is reduced to
this by consideration of nodal domains.

In Section \ref{secMFMC} we explain the Max Flow Min Cut Theorem, and
in Section \ref{secsubdomains} we show how a classical inequality bounding the Cheeger constant of a plane domain in terms of its inradius may be understood in terms of vector fields.

\section{Max Flow Min Cut Theorems} \label{secMFMC}
The classical Max Flow Min Cut Theorem deals with a {\em discrete network},
consisting of a finite set $V$ and a function $c:V\times V\to [0,\infty)$.
$c(v,w)$ may be considered as the capacity of a pipe connecting the
'nodes' $v,w\in V$.
Two nodes are distinguished,
the {\em source} $s$ and the {\em sink} $t$.

We want to transport some liquid from $s$ to $t$. A (stationary, i.e.\ time independent) 'transport plan' is modelled by a {\em flow}, i.e.\ a function $f:V\times V\to[0,\infty)$ satisfying the {\em capacity constraint} $f(v,w)\leq c(v,w)\,\forall v,w\in V$ and the 'Kirchhoff law'
\begin{equation}\label{kirchhoff}
N_f(v) := \sum_w [f(v,w)-f(w,v)] = 0\ \forall v\in V\setminus\{s,t\},
\end{equation}
that is, the total flow out of $v$ equals the total flow into $v$, except at the source and the sink.
We let $$\value(f):=N_f(s),$$
the net flow out of the source.
\eqref{kirchhoff} implies that this equals $-N_f(t)$.

The question is how big the value of a flow can be, given the capacity constraint. A simple upper bound can be given by any {\em cut}, i.e.\
subset $S\subset V$ containing $s$ but not $t$. Clearly, if we define
$$\cost(S) := \sum_{v\in S,w\in V\setminus S} c(v,w),$$
the total capacity of pipes leaving $S$, then
\begin{equation}\label{trivial}
 \value(f)\leq \cost (S)
\end{equation}
for any flow $f$ and any cut $S$: Any net flow out ouf the source must reach the sink, so it must leave $S$ at some point.

\begin{theorem}[Max Flow Min Cut Theorem, \cite{EliFeiSha:NMFTN}, \cite{ForFul:MFTN}] \label{FFMFMC}
In any network there is a flow $\fmax$ and a cut $\Smin$ satisfying
$\value(\fmax) = \cost(\Smin)$.
\end{theorem}

Various generalizations of Theorem \ref{FFMFMC} have been proposed.
First, one may allow several sources and several sinks of prescribed relative 'strengths'. As an example analogous to our continuous setup below, all but one
nodes could be sources, of equal strength, the remaining node being the sink $t$.
Then condition \eqref{kirchhoff} is empty. If we let $\value(f)$ be the minimum
net flow out of any source node then \eqref{trivial} becomes
$\value(f)\leq \cost(S)/|S|$ for any set $S$ not containing $t$, and
the corresponding Max Flow Min Cut theorem holds again.

More challenging are generalizations to infinite sets $V$, modelling continuously distributed sources or sinks, for example.
Measure theoretic versions were proved in \cite{ChaHar:DTCAFFFN} and \cite{FucLus:CC}.
%: $V$ is replaced by a measurable space $(\Omega,\Sigma), capacity and flow are %measures $\tau$,
%$\nu$ on $\Omega\times\Omega$ (with $\tau(A\times B)$, $\nu(A\times B)$ %modelling the capacity of pipes from $B$ to $A$, respectively the flow from
%$B$ to $A$), and sources and sinks are modelled by a 'consumption measure' on %$\Omega$.
Here we are more interested in a geometric model. Several, slightly different, models were proposed in \cite{Iri:TFCAFN}, \cite{Str:MFTD}, \cite{TagIri:CADNAAURN} and later unified and generalized by Nozawa \cite{Noz:MFMCTAN}. (I recommend \cite{Str:MFTD} for enjoyable reading.)

We explain the relation of network flows to
Theorem \ref{MFMCcont}, which is a special case of Nozawa's general Max Flow Min Cut Theorem:
The domain $\Omega\subset\R^n$ is the network.
A (stationary) flow is modelled by a vector field $V$, as is common in continuum fluid dynamics. The capacity constraint is \eqref{Vlength}.
Sources are distributed uniformly over $\Omega$, and \eqref{Vdiv} states that they produce liquid at a rate $h$, at least. The complement (or boundary)
of $\Omega$ should be considered as the sink (one should think of a single sink, i.e. collapse the boundary to a point; in this way one does not need to prescribe the relative strengths of the sinks along the boundary).
%alternatively, one could think of open networks, which in the discrete case %would mean taking out $t$ but leaving the edges attached to it as 'open' edges %in the network).
A cut is a subset $S\subset\Omega$, and \eqref{dSVS} states the obvious fact that anything that is produced within $S$ must leave $S$ through its boundary,
which yields the bound in Proposition \ref{prop}.

The discrete Max Flow Min Cut Theorem is proved, in most texts on discrete optimization (see \cite{KorVyg:CO}, for example), by inductively constructing $\fmax$.
However, as is already remarked in \cite{ForFul:FN}, the theorem is also an instance of the very general duality principle in convex optimization, and
this approach also yields the generalizations mentioned above.

The duality principle associates to our optimization problem
(maximize $h=\inf_{x\in\Omega}\div V(x)$, subject to the constraint \eqref{Vlength}) a dual problem, which turns out to be:
\begin{equation}\label{dual}
\text{Minimize } Q(\phi):=\frac{\|\phi\|_{BV}}{\|\phi\|_{L^1}},\
\text{ subject to } \phi\geq 0.
\end{equation}
Here, $\phi$ is a function of bounded variation on $\R^n$ which vanishes outside $\Omega$, $\|\phi\|_{BV}$ is its total variation (in $\R^n$), which equals
$\int_\Omega |\grad\phi|$ in case $\phi$ is smooth, and $\|\phi\|_{L^1}=\int_\Omega \phi$.
%
%\begin{equation}\label{dual}
%\text{Minimize \|u\|_{BV}:=\int_\Omega|\grad u|, \ \text{ subject to }
%u\geq 0, \int u = 1.
%\end{equation}
The general duality theorem says that $\sup_V \inf_x (\div V(x)) = \inf_\phi Q(\phi)$.
To see the relation to Cheeger's constant, first recall (see \cite{EvaGar:MTFPF}, Chapter 5) that the perimeter of $S\subset\R^n$ is $|\partial S| = \|\chi_S\|_{BV}$, if this is finite, where $\chi_S$ is the characteristic function of $S$. Therefore, $Q(\chi_S) = |\partial S|/|S|$
if $S\subset\Omega$.
Next, the coarea formula (loc.\ cit.) states that
for $\phi\geq 0$ of bounded variation, the sets $\{\phi>t\}$ have finite perimeter for almost all $t$, and
\newcommand{\chit}{\chi_{\{\phi>t\}}}
\begin{equation}\label{coarea}
\|\phi\|_{BV} = \int_0^\infty \|\chit\|_{BV}\,dt.
\end{equation}
Since $\phi=\int_0^\infty \chit\,dt$, one also has $\|\phi\|_{L^1}=\int_0^\infty
\|\chit\|_{L^1}\,dt$, and therefore $Q(\phi)\geq \inf_t Q(\chit)$.
This shows that in \eqref{dual} one may restrict $\phi$ to characteristic functions, so the infimum is precisely Cheeger's constant.

\section{Cheeger's constant and inradius}\label{secsubdomains}
In this section let $\Omega$ be a simply connected plane domain, and
let $\rho_\Omega$ denote its inradius. Also, define the 'reduced inradius'
\begin{equation}\label{rhotilde}
\rhotilde_\Omega := \frac{\rho_\Omega}{1+\pi\rho_\Omega^2/|\Omega|}.
\end{equation}
Clearly, $\rho_\Omega/2 < \rhotilde_\Omega < \rho_\Omega$.
A well-known lower bound for $\lambda_\Omega$ is
\begin{equation}\label{Ossermanineq}
\lambda_\Omega\geq \frac1{4\rhotilde_\Omega^2}.
\end{equation}
The weaker estimate $\lambda_\Omega>1/4\rho_\Omega^2$ is sometimes called Hayman's inequality or Osserman's inequality, since
it was proved by Hayman (with 4 replaced by 900) \cite{Hay:SBPF} and Osserman  \cite{Oss:NHTBND}, but in fact it was first proved by E.\ Makai \cite{Mak:LEPFSCM}.\footnote{The constant 1/4 in $\lambda_\Omega>1/4\rho_\Omega^2$ has since been improved using ideas from probability and conformal mapping, the currently best value is
0.6197, see \cite{BanCar:BMFFD}.}
There are similar estimates in the multiply connected case (to which the considerations below apply as well), but there is no direct higher dimensional generalization. See \cite{Bog:IISNEP} for a generalization to a certain pseudo-Laplacian.

\eqref{Ossermanineq} follows from
\begin{equation}\label{OsChineq}
h_\Omega\geq \frac1\rhotilde_\Omega
\end{equation}
and Cheeger's inequality. Note that \eqref{OsChineq} is sharp for the disk.
\eqref{OsChineq} is implicit in \cite{Mak:LEPFSCM} and \cite{Oss:NHTBND}, but does not seem to be stated explicitly in the literature.
Let us give the proof along the lines of \cite{Mak:LEPFSCM}, \cite{Oss:NHTBND}:
Let
\begin{equation}\label{classS}
 \calS=\{S\subset\Omega:\, S\text{ open and simply connected},\, \partial S\text{ smooth}\}.
\end{equation}
Clearly,
\begin{equation}\label{subdomcond}
 h_{\Omega} = \inf_{S\in\calS} \frac {|\partial S|} {|S|}
\end{equation}
since filling in all 'holes' in an arbitrary $S\subset \Omega$, making it simply connected, increases $|S|$, decreases $|\partial S|$, and results in a subset of $\Omega$ (since $\Omega$ is simply connected).
Also,
\begin{equation}\label{monotone}
S\subset\Omega \Rightarrow \rhotilde_S\leq\rhotilde_\Omega.
\end{equation}
To see this, first note that, for $A,\rho>0$, the function $f_A(\rho)=\rho\mapsto \rho/(1+\pi\rho^2/A)$ is increasing in $\rho$ for $\pi\rho^2\leq A$.
Now $|S|\leq|\Omega|$ yields $\rhotilde_S=f_{|S|}(\rho_S)\leq f_{|\Omega|}(\rho_S)$; also $\pi\rho_S^2\leq \pi\rho_\Omega\leq|\Omega|$,
so the monotonicity of $f_{|\Omega|}$ gives \eqref{monotone}.
Now the main step is 'Bonnesen's inequality': For simply connected $S\subset\R^2$
\begin{equation}{\label{bonn}}
\rho_S|\partial S|\geq |S|+\pi\rho_S^2.
\end{equation}
This, together with \eqref{subdomcond} and \eqref{monotone} proves
\eqref{OsChineq}.\footnote{Note that \eqref{bonn} is equivalent to $|\partial S|^2 - 4\pi|S| \geq (|\partial S|-2\pi\rho_S)^2$ and therefore a sharper version of the classical isoperimetric inequality $|\partial S|^2\geq 4\pi|S|$. \eqref{bonn}
was proved by Bonnesen for convex $\Omega$
\cite{Bon:UVIUKEKAMUKK} and by Besicovitch \cite{Bes:VCIP} for general simply connected domains, see also Sz.-Nagy \cite{SzNag:UPNEB}.
}
%Idea: If the boundary of the  $t$-level-set is smooth then
%$ \frac{d^2}{dt^2} A_t = -\int_{T_t}\kappa = -2\pi k < 0$ where $\kappa$ is the curvature and $k$ the number of components; in general, this boundary is not smooth (not even for almost all $t$), so a suitable approximation argument is used.)}
\medskip

The question arises naturally whether one may prove \eqref{OsChineq} by constructing a vector field $V$ on $\Omega$ satisfying \eqref{Vlength}, \eqref{Vdiv} with $h=1/\rhotilde_\Omega$. It is not clear how to do this. It seems more interesting to infer from Theorem \ref{MFMCcont} and \eqref{OsChineq}:

\begin{corollary}\label{Osscor}
Let $\Omega$ be a simply connected plane domain with Lipschitz boundary and reduced inradius $\rhotilde_\Omega$ defined by \eqref{rhotilde}. Then
there is a vector field $V$ on $\Omega$ satisfying \eqref{Vlength}, \eqref{Vdiv}
with $h=1/\rhotilde_\Omega$.
\end{corollary}

Although there seems to be no natural, geometrically defined candidate for this vector field, we now proceed to show how  certain geometric vector fields for {\em subdomains} of $\Omega$ yield \eqref{OsChineq}.

First, we need the following variant of Proposition \ref{prop}:

\begin{proposition} \label{subdomains}
Let $\Omega\subset\R^n$ be open, and let $\calS$ be a class of Lipschitz subdomains of
$\Omega$ satisfying \eqref{subdomcond}.
Let $h\in\R$.
Suppose that for each $S\in\calS$ there is a vector field $V_S\in \XBV(S)$ on $S$ satisfying
\begin{eqnarray}
|V_S| &\leq & 1,\quad\text{pointwise on }|S|,  \label{VSlength}\\
\int_S \div V_S&\geq & h|S| \label{VSdiv}.
\end{eqnarray}
Then $h_\Omega \geq h$.
\end{proposition}
\proof
This follows from \eqref{dSVS} applied to $V_S$.
\endproof

Note that condition \eqref{VSdiv} is weaker than the pointwise condition \eqref{Vdiv} (if this is applied to $V_S$ on $S$). So in order to get effective lower bounds on $h_\Omega$, one has more flexibility in choosing 'test' vector fields, but one needs to do it for all $S\in\calS$ simultaneously.

As before, we define $\calS$ by \eqref{classS}. Using \eqref{monotone}, we then see that
\eqref{OsChineq} follows from Proposition \ref{subdomains} and the following:

\begin{proposition}\label{propinr}
Let $S\subset\R^2$ be a smooth, simply connected domain, of inradius $\rho_S$ and reduced inradius $\rhotilde_S$. Then
there is a vector field $V\in\XBV(S)$ on $S$ satisfying
\begin{eqnarray}
|V| &\leq & 1,\\
\int_S \div V&\geq & \frac{|S|}{\rhotilde_S} = \frac{|S|}{\rho_S} + \pi\rho_S. \label{VSdiv2}
\end{eqnarray}
\end{proposition}
\proof
By scale invariance we may assume $\rho_S=1$.
Let $\phi(x)$ denote the distance of $x\in S$ to $\partial S$, and define
\begin{equation}\label{vect}
 V=\frac12 \grad (1-\phi)^2 = -(1-\phi)\grad \phi.
\end{equation}
This is motivated by the case of the disk, Example \ref{ex}.

We will use the following facts about the distance function:
\begin{enumerate}
\item[(a)] $\grad\phi\in L^\infty(S,\R^n)$, and $|\grad\phi|=1$ almost everywhere.
\item[(b)] $\grad \phi$ has bounded variation.
\item[(c)] For almost all $t\in[0,1]$ the level set $\phit$ consists of a finite union of piecewise smooth, simple closed curves, with non-zero angles, and
$\partial\{\phi<t\} = \partial S\cup \phit.$
\end{enumerate}
(a) and (c) are proved, for example, in \cite{Har:GPCL}, by a detailed analysis of the function $F:[0,L]\times [0,1]\to \R^2$, with $L$ the length of $\partial S$, defined by requiring that $s\mapsto F(s,0)$ is an arclength parametrization of $\partial S$ and $t\mapsto F(s,t)$ is the unit speed normal to $\partial S$ starting inward at $F(s,0)$. For higher dimensional generalizations of (c), stating that $\phit$ is Lipschitz for a.e.\ $t$ and using the notion of Clarke gradient of $\phi$, see \cite{ItoTan:STDF}, \cite{Rif:MSTDFRM}.
(b) is folklore\footnote{Proof: Let $x_0\in S$ and $a=\phi(x_0)>0$.
By an easy calculation, $x\mapsto \phi_y(x):= |x-y|-\frac2a |x-x_0|^2$ has negative definite Hessian at $x_0$ (if $y\in\partial S$) and so is concave near $x_0$. Therefore, $\tilde{\phi}(x)=\inf_{y\in\partial S}\phi_y(x)$ is concave near $x_0$, so $\grad\tilde{\phi}$ has bounded variation near $x_0$
(see \cite{EvaGar:MTFPF}, Section 6.3, Theorem 3).
Finally, $\phi(x)=\tilde{\phi}(x)+\frac2a |x-x_0|^2$ shows that $\grad\phi$ has locally bounded variation in $S$, and since $\phi$ is smooth near $\partial S$ it has bounded variation.}.

Therefore, to prove Proposition \ref{propinr}
 it remains to verify $\int\div V \geq |S|+\pi$.
Now $\div V= |\grad \phi|^2 - (1-\phi)\Delta \phi$, with $\Delta$ the Laplace operator. Since $|\grad\phi|=1$ we need to show
\begin{equation}\label{finalest}
\int_S (1-\phi)\Delta \phi \leq -\pi.
\end{equation}
Now we have $\int_S (1-\phi)\Delta \phi
=\int_S(\int_{\phi(x)}^1dt)\Delta\phi
=\int_0^1 (\int_{\philt}\Delta\phi)\,dt$, by applying
Fubini's theorem for measures (see \cite{EvaGar:MTFPF}, for example)
to $\int_U \Delta\phi\,dt$, where $U=\{(t,x):\phi(x)<t\}\subset [0,1]\times S$.
%Let $S_t=\{x\in S:\, d(x)<t\}$. Then $S_t=S$ for $t\geq 1$ (since the inradius
%is one)
%$$\int_S (1-d)\Delta d = \int_0^1\left(\int_{S_t} \Delta d\right)\, dt.$$
If $t$ is as in (c) above then,
by the divergence theorem, $\int_{\phi<t}\Delta \phi = \int_{\partial\{\phi<t\}}
(\grad \phi)\cdot dn=L_t-L_0$, where $L_t$ denotes the length of $\phit$.
Finally, $\int_0^1L_t\,dt = |S|$ (by the coarea formula \eqref{coarea}, for example, using (a)), and therefore
$$\int_0^1 (\int_{\philt}\Delta\phi)\,dt = \int_0^1 (L_t-L_0)\,dt=|S|-|\partial S|\leq-\pi$$
by \eqref{bonn}, and this proves \eqref{finalest}.
\endproof

Note that the vector field \eqref{vect} cannot be used directly in Proposition
\ref{prop}, since $\Delta\phi>0$ near concave parts of $\partial\Omega$, so that the pointwise estimate $\div V\geq1$ is false.

{\bf Remark:}
If $\phi$ was smooth everywhere, one could prove \eqref{finalest} without appealing to Bonnesen's inequality (and, in effect, reprove this inequality), as follows: The general coarea formula gives
$$\int_S (1-\phi)\Delta \phi = \int_0^1 (\int_\phit \frac{1-\phi}{|\grad\phi|}\Delta\phi)\,dt = \int_0^1(1-t) (\int_\phit \Delta\phi)\,dt,$$
where the line integrals are with respect to arclength measure.
It is elementary to see that $\Delta\phi$ equals
minus the curvature of $\phit$ (whereever $\phi$ is smooth). Also, the integral of the curvature along a smooth simple closed curve equals $2\pi$. So if almost all level sets were smooth (instead of piecewise smooth) we would obtain
\eqref{finalest}.

The problem with this 'proof' is that, typically,  $\Delta\phi$ is not a function but a measure (so that the coarea formula is not applicable),
and a positive measure of level sets may be non-smooth.
Consider, for example, a rectangle: $\grad\phi$ has a jump at its 'center line', leading to a $\delta$ type singularity of $\Delta\phi$ there.

It should be possible and would be interesting to find generalizations of the coarea formula
(using suitable transversality hypotheses) and of the curvature argument
that make this proof work.

\bibliography{dglib}
% bbl-File wurde von Hand geaendert (review lines auskomnentiert)
\end{document}